\documentclass[11pt]{amsart}
\usepackage{amssymb}
\usepackage{amsmath}
\usepackage{fancyhdr}
\usepackage[british]{babel}
\usepackage{geometry}
\usepackage{enumitem}
\usepackage{algpseudocode}
\usepackage{dsfont}
\usepackage{centernot}
\usepackage{xstring}
\usepackage{colortbl}
\usepackage[section]{algorithm}
\usepackage{graphicx}
\usepackage[font={footnotesize}]{caption}
\usepackage[usenames,dvipsnames,table]{xcolor}
\usepackage{pstricks,tikz}
\usetikzlibrary{arrows,decorations.pathmorphing}
\usepackage[h]{esvect}
\usepackage[normalem]{ulem}
\usepackage[
		bookmarksopen=true,
		bookmarksopenlevel=1,
		colorlinks=true,
		linkcolor=darkblue,
        linktoc=page,
		citecolor=darkblue,
]{hyperref}

\title[]{Minimalist designs}
\date{\today}
\author[B. Barber, S. Glock, D. K\"uhn, A. Lo, R. Montgomery and D. Osthus]{Ben Barber, Stefan Glock, Daniela K\"uhn, Allan Lo, Richard Montgomery and Deryk Osthus}

\thanks{The research leading to these results was partially supported by the European Research Council under the European Union's Seventh Framework Programme (FP/2007--2013) / ERC Grant 306349 (S.~Glock and D.~Osthus). The research was also partially supported by the EPSRC, grant nos. EP/N019504/1 (D.~K\"uhn) and EP/P002420/1 (A.~Lo), and by the Royal Society and the Wolfson Foundation (D.~K\"uhn).}

\newtheorem{theorem}[algorithm]{Theorem}

\newtheorem{lemma}[algorithm]{Lemma}
\newtheorem{cor}[algorithm]{Corollary}
\newtheorem{fact}[algorithm]{Fact}

\theoremstyle{definition}

\newtheorem{conj}[algorithm]{Conjecture}
\newtheorem{defin}[algorithm]{Definition}

\newtheoremstyle{claimstyle}{5pt}{5pt}{\em}{5pt}{\em}{:}{5pt}{}
\theoremstyle{claimstyle}

\newtheoremstyle{stepstyle}{10pt}{5pt}{\em}{0pt}{\em}{:}{5pt}{}
\theoremstyle{stepstyle}

\numberwithin{equation}{section}

\definecolor{darkblue}{rgb}{0,0,0.5}

\def\noproof{{\unskip\nobreak\hfill\penalty50\hskip2em\hbox{}\nobreak\hfill%
       $\square$\parfillskip=0pt\finalhyphendemerits=0\par}\goodbreak}
\def\endproof{\noproof\bigskip}

\newdimen\margin
\def\textno#1&#2\par{
   \margin=\hsize
   \advance\margin by -4\parindent
          \setbox1=\hbox{\sl#1}
   \ifdim\wd1 < \margin
      $$\box1\eqno#2$$
   \else
      \bigbreak
      \hbox to \hsize{\indent$\vcenter{\advance\hsize by -3\parindent
      \it\noindent#1}\hfil#2$}
      \bigbreak
   \fi}

\def\proof{\removelastskip\penalty55\medskip\noindent\setcounter{claim}{0}\setcounter{step}{0}{\bf Proof. }} 

\def\lateproof#1{\removelastskip\penalty55\medskip\noindent\setcounter{claim}{0}\setcounter{step}{0}{\bf Proof of #1. }} 

\begin{document}

\newcommand{\new}[1]{\textcolor{red}{#1}}
\def\COMMENT#1{}
\def\TASK#1{}

\newcommand{\todo}[1]{\begin{center}\textbf{to do:} #1 \end{center}}

\def\eps{{\varepsilon}}
\newcommand{\ex}{\mathbb{E}}
\newcommand{\pr}{\mathbb{P}}
\newcommand{\cB}{\mathcal{B}}
\newcommand{\cA}{\mathcal{A}}
\newcommand{\cE}{\mathcal{E}}
\newcommand{\cS}{\mathcal{S}}
\newcommand{\cF}{\mathcal{F}}
\newcommand{\cG}{\mathcal{G}}
\newcommand{\bL}{\mathbb{L}}
\newcommand{\bF}{\mathbb{F}}
\newcommand{\bZ}{\mathbb{Z}}
\newcommand{\cH}{\mathcal{H}}
\newcommand{\cC}{\mathcal{C}}
\newcommand{\cM}{\mathcal{M}}
\newcommand{\bN}{\mathbb{N}}
\newcommand{\bR}{\mathbb{R}}
\def\O{\mathcal{O}}
\newcommand{\cP}{\mathcal{P}}
\newcommand{\cQ}{\mathcal{Q}}
\newcommand{\cR}{\mathcal{R}}
\newcommand{\cJ}{\mathcal{J}}
\newcommand{\cL}{\mathcal{L}}
\newcommand{\cK}{\mathcal{K}}
\newcommand{\cD}{\mathcal{D}}
\newcommand{\cI}{\mathcal{I}}
\newcommand{\cV}{\mathcal{V}}
\newcommand{\cT}{\mathcal{T}}
\newcommand{\cU}{\mathcal{U}}
\newcommand{\cX}{\mathcal{X}}
\newcommand{\cZ}{\mathcal{Z}}
\newcommand{\1}{{\bf 1}_{n\not\equiv \delta}}
\newcommand{\eul}{{\rm e}}
\newcommand{\Erd}{Erd\H{o}s}
\newcommand{\cupdot}{\mathbin{\mathaccent\cdot\cup}}
\newcommand{\whp}{whp }
\newcommand{\bX}{\mathcal{X}}
\newcommand{\bV}{\mathcal{V}}

\newcommand{\doublesquig}{%
  \mathrel{%
    \vcenter{\offinterlineskip
      \ialign{##\cr$\rightsquigarrow$\cr\noalign{\kern-1.5pt}$\rightsquigarrow$\cr}%
    }%
  }%
}

\newcommand{\defn}{\emph}

\newcommand\restrict[1]{\raisebox{-.5ex}{$|$}_{#1}}

\newcommand{\prob}[1]{\mathrm{\mathbb{P}}\left[#1\right]}
\newcommand{\expn}[1]{\mathrm{\mathbb{E}}\left[#1\right]}
\def\gnp{G_{n,p}}
\def\G{\mathcal{G}}
\def\lflr{\left\lfloor}
\def\rflr{\right\rfloor}
\def\lcl{\left\lceil}
\def\rcl{\right\rceil}

\newcommand{\qbinom}[2]{\binom{#1}{#2}_{\!q}}
\newcommand{\binomdim}[2]{\binom{#1}{#2}_{\!\dim}}

\newcommand{\grass}{\mathrm{Gr}}

\newcommand{\brackets}[1]{\left(#1\right)}
\def\sm{\setminus}
\newcommand{\Set}[1]{\{#1\}}
\newcommand{\set}[2]{\{#1\,:\;#2\}}
\newcommand{\krq}[2]{K^{(#1)}_{#2}}
\newcommand{\ind}[1]{$\mathbf{S}(#1)$}
\newcommand{\indcov}[1]{$(\#)_{#1}$}
\def\In{\subseteq}

\begin{abstract}  \noindent
The iterative absorption method has recently led to major progress in the area of (hyper-)graph decompositions. Amongst other results, a new proof of the Existence conjecture for combinatorial designs, and some generalizations, was obtained. Here, we illustrate the method by investigating triangle decompositions:
we give a simple proof that a triangle-divisible graph of large minimum degree has a triangle decomposition and prove a similar result for quasi-random host graphs.
\end{abstract}

\maketitle

\section{Introduction} \label{sec:intro}

\subsection{Steiner triple systems and the decomposition threshold}
A famous theorem of Kirkman~\cite{kirkman:47} from 1847 states that a Steiner triple system of order~$n$ exists if and only if $n\equiv 1,3\mod{6}$. Here, a \defn{Steiner triple system of order~$n$} is a collection of $3$-subsets of $[n]$ such that every $2$-subset of $[n]$ is contained in exactly one of the $3$-sets. More generally, given graphs $G$ and $F$, we say that $G$ has an \defn{$F$-decomposition} if its edge set can be partitioned into copies of $F$. Clearly, a Steiner triple system of order~$n$ is equivalent to a $K_3$-decomposition of~$K_n$.
Observe that if a graph $G$ admits a $K_3$-decomposition, then the number of edges of $G$ must be divisible by~$3$, and all the vertex degrees of $G$ must be even. We say that $G$ is \defn{$K_3$-divisible} if it has these properties. Clearly, not every $K_3$-divisible graph has a $K_3$-decomposition (e.g.~$C_6$). In fact, to decide whether a given graph has a $K_3$-decomposition is NP-hard~\cite{DT:97}.
However, the following beautiful conjecture of Nash-Williams suggests that if the minimum degree of $G$ is sufficiently large, then divisibility is not only necessary but also sufficient for the existence of a $K_3$-decomposition.

\begin{conj}[Nash-Williams~\cite{nash-williams:70}]  \label{conj:NW}
\emph{For sufficiently large $n$, every $K_3$-divisible graph $G$ on $n$ vertices with $\delta(G)\ge 3n/4$ has a $K_3$-decomposition.}
\end{conj}

The following class of extremal examples shows that the bound on the minimum degree would be best possible. Given any $k\in \bN$, let $G_1$ and $G_2$ be vertex-disjoint $(6k+2)$-regular graphs with $|G_1|=|G_2|=12k+6$ and let $G_3$ be the complete bipartite graph between $V(G_1)$ and $V(G_2)$. Let $G:=G_1\cup G_2\cup G_3$. (In the standard construction, each of $G_1$ and $G_2$ is a union of two disjoint cliques of size~$6k+3$.)
Clearly, $\delta(G)=3|G|/4-1$\COMMENT{$|G|=24k+12$. $\delta(G)=6k+2+12k+6=18k+8$.} and $G$ is $K_3$-divisible.\COMMENT{all degrees even, $e(G_3)=(12k+6)^2$, $e(G_1\cup G_2)=(12k+6)(6k+2)$} However, every triangle in $G$ contains at least one edge from $G_1\cup G_2$. Since $2e(G_1\cup G_2)<e(G_3)$, $G$ cannot have a $K_3$-decomposition.

For $n\in \bN$, define $\delta(n)$ as the minimum of all natural numbers $d\in \bN$ such that every $K_3$-divisible graph $G$ on $n$ vertices with $\delta(G)\ge d$ has a $K_3$-decomposition. Hence, Conjecture~\ref{conj:NW} is equivalent to saying that $\delta(n)\le 3n/4$ for all sufficiently large~$n$.

The \defn{decomposition threshold of $K_3$} is defined as $\delta_{K_3}:= \limsup_{n\to \infty} \frac{\delta(n)}{n}$. Conjecture~\ref{conj:NW} would imply that $\delta_{K_3}\le 3/4$ (and equality would hold by the above example).
More generally, the \defn{$F$-decomposition threshold} $\delta_F$ can be defined analogously for any \mbox{(hyper-)}graph~$F$ we wish to decompose into.
It is conjectured in~\cite{GKLMO:19} that for every graph~$F$, we have $\delta_F\le 1-1/(\chi(F)+1)$ (but equality does not hold for every~$F$).
This is reminiscent of related results in Extremal Combinatorics. For example, the Erd\H{o}s--Stone--Simonovits theorem says that any graph $G$ with $e(G)\ge (1-1/(\chi(F)-1)+o(1))\binom{n}{2}$ contains $F$ as a subgraph, and a theorem due to Alon and Yuster states that if $\delta(G)\ge (1-1/\chi(F)+o(1))n$ and $|F|$ divides $|G|$ then $G$ has an $F$-factor (again, the latter bound is not optimal for every~$F$, see~\cite{KO:09}).

The main application of the iterative absorption method is to turn an approximate decomposition into an exact decomposition. To formalize this, we define the approximate
decomposition threshold  $\delta^{0+}$ to be the infimum of all $\delta\in[0,1]$ with the following property: for all $\gamma>0$, there exists $n_0\in \bN$ such that any graph $G$ on $n\ge n_0$ vertices with $\delta(G)\ge \delta n$ contains a $K_3$-decomposable subgraph~$H$ such that $\Delta(G-H)\le \gamma n$.
In this paper, we prove the following theorem, which reduces the challenge of showing that $\delta_{K_3}=3/4$ to showing that $\delta^{0+}\le 3/4$. The result itself is already contained in~\cite{BKLO:16}. However, the proof we present here is much simpler and serves as an illustration as to how to use iterative absorption for decomposition problems.

\begin{theorem} \label{thm:main}
Let $\delta:=\max\Set{3/4,\delta^{0+}}$.
For all $\eps>0$, there exists $n_0\in \bN$ such that every $K_3$-divisible graph $G$ on $n \ge n_0$~vertices with $\delta(G)\ge (\delta+\eps)n$ has a $K_3$-decomposition.
\end{theorem}

In a nutshell, Theorem~\ref{thm:main} says that $\delta_{K_3}\le \max\Set{3/4,\delta^{0+}}$.
Any improvement on the value of $\delta^{0+}$ thus immediately improves the value of~$\delta_{K_3}$. It follows from results in~\cite{DP:19,HR:01} that $\delta^{0+}\le (7+\sqrt{21})/14$. Together with Theorem~\ref{thm:main} this implies that $\delta_{K_3}\le (7+\sqrt{21})/14\approx 0.82733$. It would be very desirable to prove that $\delta^{0+}\le 3/4$, which would thus give an asymptotic version of Conjecture~\ref{conj:NW}.
(We will discuss $\delta^{0+}$ further in Section~\ref{sec:fractional}.)
In order to prove Conjecture~\ref{conj:NW} in full, however, would in addition likely require a stability analysis --- a daunting prospect given the abundance of extremal examples. On the other hand, for decompositions into even cycles (except of length~$6$), this has been carried out in~\cite{taylor:19}.
Further results in the spirit of Theorem~\ref{thm:main} for decompositions into arbitrary graphs~$F$ were obtained in~\cite{BKLO:16,GKLMO:19}. In particular, in~\cite{GKLMO:19} $\delta_F$ is determined for every bipartite graph~$F$, and it is shown that the threshold of cliques equals its fractional version.

\subsection{Fractional decompositions} \label{sec:fractional}
A successful approach to obtain bounds on $\delta^{0+}$ is to take a detour via fractional decompositions. Let $\binom{G}{F}$ be the set of copies of~$F$ in $G$.
A \defn{fractional $F$-decomposition of $G$} is a function $\omega \colon \binom{G}{F} \to [0,1]$ such that for all $e \in E(G)$,
\begin{equation} \label{packing}
\sum_{F' \in \binom{G}{F} \colon e \in E(F')} \omega(F') = 1.
\end{equation}
Thus, an $F$-decomposition is a fractional $F$-decomposition with image~$\Set{0,1}$.

We define the \emph{fractional decomposition threshold $\delta^\ast$ of $K_3$} to be the infimum of all $\delta\in[0,1]$ with the following property: there exists $n_0\in \bN$ such that any graph $G$ on $n\ge n_0$ vertices with $\delta(G)\ge \delta n$ has a fractional $K_3$-decomposition.
Observe that the extremal example after Conjecture~\ref{conj:NW} also implies that $\delta^\ast\ge 3/4$.

Haxell and R\"odl~\cite{HR:01} showed  that if $G$ is an $n$-vertex graph with a fractional $F$-decomposition, then all but $o(n^2)$ edges of $G$ can be covered with edge-disjoint copies of~$F$.
Note that the definition of $\delta^{0+}$ requires a leftover of small maximum degree, whereas this result only provides a leftover with $o(n^2)$ edges. It is however easy to turn such a leftover into one with small maximum degree
(see e.g.~Lemma~10.6 in~\cite{BKLO:16}).
Thus, $\delta^{0+}\le \delta^\ast$.

Until recently, the best bound on $\delta^\ast$ was obtained by Dross~\cite{dross:16},
who showed that $\delta^\ast \le 0.9$, using an elegant approach based on the
max-flow-min-cut theorem. 
Very recently, Delcourt and Postle~\cite{DP:19} showed that $\delta^\ast \le (7+\sqrt{21})/14 \approx 0.82733$. (A slightly weaker bound was obtained simultaneously by Dukes and Horsley~\cite{DH:19}.)
The best current bound on the
fractional decomposition threshold of larger cliques was proved in~\cite{montgomery:17}.
It would be very interesting to improve these results.

\subsection{Quasi-random host graphs}
Rather than graphs of large minimum degree, it also makes sense to consider
quasi-random host graphs.
A natural notion of quasi-randomness in this context is that of typicality.
Given $p,\xi>0$ and $h \in\mathbb{N}$,  an $n$-vertex graph $G$ is
\emph{$(\xi,h,p)$-typical} if
for every set $A \subseteq V(G)$ with $|A| \le h$ the common neighbourhood
of the vertices in $A$ has size $(1\pm \xi)p^{|A|}n$.
Note that a binomial random graph with edge probability $p$ is likely to be typical.
The following result
(with $h$ large) was first proved by Keevash~\cite{keevash:14}.
\begin{theorem} \label{thm:qr}
For all $p>0$, there exist $n_0\in \bN$ and $\xi>0$ such that every $(\xi,4,p)$-typical $K_3$-divisible graph on $n \ge n_0$~vertices has a $K_3$-decomposition.
\end{theorem}
In Section~\ref{sec:qr}, we will outline how the proof of Theorem~\ref{thm:main}
can be adapted to give a proof of Theorem~\ref{thm:qr}.

\subsection{Designs} \label{designs}
More generally, it is natural to consider these questions for uniform hypergraphs.
In particular, such decompositions give rise to combinatorial designs with arbitrary parameters. More precisely, let $F,H$ be $r$-graphs (i.e.~$r$-uniform hypergraphs) and let $K^{r}_{n}$ denote the complete $r$-graph on $n$ vertices.
An \defn{$F$-decomposition} of $H$ is a collection  of copies of
$F$ in $H$ such that every edge of $H$ is contained in exactly one of these copies. A \defn{$(q,r,\lambda)$-design} of $H$ is a collection  of distinct copies of $K^{r}_{q}$ in $H$ such that every edge of $H$ is contained in exactly $\lambda$ of these copies. (When $\lambda=1$ and $H=K^{r}_{n}$, these are referred to as \emph{Steiner systems}.)

The existence of  $(q,r,\lambda)$-designs of $K^{r}_{n}$ for fixed $q,r,\lambda$
was proved by Keevash~\cite{keevash:14} (subject to $n$ satisfying the necessary divisibility conditions).
Much more generally, he proved the existence of $(q,r,\lambda)$-designs of $H$,
where $H$ is quasi-random (and dense).
His proof is based on algebraic and probabilistic techniques
(see~\cite{keevash:18} for an exposition of the triangle case and~\cite{keevash:18c} for further discussions).

A new proof (based on iterative absorption) of this result was given in~\cite{GKLO:16}.
More generally,~\cite{GKLO:16} proves the existence of $F$-designs for arbitrary $r$-graphs $F$ and 
also provides bounds on the decomposition threshold of~$K^{r}_{q}$ as well as more general $r$-graphs $F$. 
Further results (including the existence of resolvable designs and a new proof of the existence of $F$-designs)
were subsequently proved by Keevash~\cite{keevash:18b}.

\subsection{Iterative absorption}
The main idea of the absorbing technique is quite natural: suppose we want to find some spanning structure in a graph or hypergraph, for instance a perfect matching, a Hamilton cycle, or an $F$-factor. In many such cases, it is much easier to find an `almost-spanning' structure, e.g.~a matching which covers almost all the vertices.
The idea of the absorbing technique is that before finding the almost-spanning structure we set aside an absorbing structure which is capable of `absorbing' the leftover vertices into the almost-spanning structure to obtain the desired spanning structure.
Such an approach was introduced systematically in the influential paper by R\"odl, Ruci\'nski and Szemer\'edi~\cite{RRS:06} to prove an analogue of Dirac's theorem for $3$-graphs (but goes back further than this, see e.g.~the work of Krivelevich~\cite{krivelevich:97} on triangle factors in random graphs, and the result of Erd\H{o}s, Gy\'arf\'as and Pyber~\cite{EGP:91} on vertex coverings with monochromatic cycles). Since then, the absorbing technique has been successfully applied to a wealth of problems concerning spanning structures.
Of course, the success of the approach stands and falls with the ability to find this `magic' absorbing structure. One key factor in this is the number of possible leftover configurations. Intuitively, the more possible leftover configurations there are, the more difficult it is to find an absorbing structure which can deal with all of them. Loosely speaking, this makes it much harder (if not impossible) to directly apply an absorbing technique for edge-decomposition problems (see e.g.~\cite[p.~343]{BKLO:16} for a back-of-the-envelope calculation).

The `iterative absorption' method overcomes this issue by splitting up the absorbing process into many steps. In each step, the number of possible leftover configurations is drastically reduced using a `partial absorbing procedure', until eventually one has sufficient control over the leftover to absorb it completely in a final absorption step. This approach was first used in~\cite{KO:13} to find Hamilton decompositions of regular robust expanders. An iterative procedure using partial absorbers was also used in \cite{KKO:15} to find optimal Hamilton packings in random graphs (though strictly speaking this is not a decomposition result).
In the context of $F$-decompositions, the method was first applied in~\cite{BKLO:16} to find $F$-decompositions of graphs of suitably high minimum degree. The results from \cite{BKLO:16} were strengthened in \cite{GKLMO:19}. Even though the overall proof in \cite{GKLMO:19} is more involved, the iterative absorption procedure itself is simpler than in~\cite{BKLO:16}.
The method has also been successfully applied to verify the Gy\'arf\'as-Lehel tree packing conjecture for bounded degree trees~\cite{JKKO:ta},
as well as to find decompositions of dense graphs in the partite setting~\cite{BKLOT:17}.
Last but not least, as mentioned in Section~\ref{designs}, the method was developed for hypergraphs (and thus~designs) in~\cite{GKLO:16}.

\subsection{Overview of the argument}

At the beginning of the proof we will fix a suitable nested sequence of vertex sets
$V(G)=U_0 \supseteq U_1 \supseteq \dots \supseteq U_\ell$, which will be called a `vortex' in $G$.
We then remove an `absorber' $A$ from $G$, described in more detail below.
After the $i$th step of our proof we can ensure that the remaining uncovered edges all lie in $U_i$, which is much smaller than $U_{i-1}$. We achieve this `cover down' step by first finding an approximate decomposition of the current leftover, and then carefully covering any 
remaining edges which do not lie inside $U_{i}$ (see Lemma~\ref{lem:cover-down}).
We can also preserve the relative minimum degree of the leftover graph $G_i$ after the $i$th iteration, 
i.e.~$\delta(G_i[U_i]) \ge (\delta+\eps/2)|U_i|$, which enables us to repeat the iteration.

The final set $U_\ell$ in the iteration will have bounded size.
This immediately implies that there are only a bounded number of possibilities $L_1,\dots,L_s$ for $L$.
We will construct the absorber $A$ for $L$ as an edge-disjoint union of 
`exclusive' absorbers $A_1, \dots, A_s$, where each $A_i$ can absorb $L_i$; that is, both $A_i$ and $A_i \cup L_i$ have a triangle decomposition for each $i \in [s]$.
Then~$A \cup L_i$ clearly has a triangle decomposition for any of the 
permissible leftovers $L_i$. 
 
 We will construct the (exclusive) absorbers in Section~\ref{absorbers}: 
rather than constructing $A_i$ directly, we will obtain it as the concatenation
of several `transformers' $T$.
The role of $T$ is  to transform $L_i$ into a suitable different graph $L_i'$
(more precisely, both $L_i' \cup T$ and $T \cup L_i$
have a triangle decomposition).
We can then concatenate several such transformers (by taking their edge-disjoint union)  to transform $L_i$ into a disjoint union of triangles,
which is clearly decomposable.


\section{Preliminaries} \label{sec:tools}
For a graph $G$, we let $|G|$ denote the number of vertices of $G$ and $e(G)$ the number of edges of $G$. We will sometimes identify a graph with its edge set if this enhances readability and does not affect the argument. For a vertex $v\in V(G)$, we write $N_G(v)$ for the neighbourhood of $v$ and $d_G(v)$ for its degree. More generally, for a subset $X\In V(G)$, we let $d_G(x,X)$ denote the number of neighbours of $x$ in~$X$.
Let $G$ be a graph and let $X,Y$ be disjoint subsets of~$V(G)$. We write $G[X]$ for the subgraph of $G$ induced by $X$, and $G[X,Y]$ for the bipartite subgraph of $G$ induced by $X,Y$.
If $G$ is a graph and $H$ is a subgraph of~$G$, then $G-H$ denotes the graph with vertex set $V(G)$ and edge set $E(G)\sm E(H)$.

We write $[n]$ for the set $\Set{1,\dots,n}$.
The expression $a= b \pm c$ means that $a\in [b-c,b+c]$.
We write $x\ll y$ to mean that for any $y\in (0,1]$ there exists an $x_0\in (0,1)$ such that for all $x\le x_0$ the subsequent statement holds. Hierarchies with more constants are defined in a similar way and are to be read from the right to the left. We will always assume that the constants in our hierarchies are reals in $(0,1]$. Moreover, if $1/x$ appears in a hierarchy, this implicitly means that $x$ is a natural number. More precisely, $1/x\ll y$ means that for any $y\in (0,1]$ there exists an $x_0\in \bN$ such that for all $x\in \bN$ with $x\ge x_0$ the subsequent statement holds.

Let $m,n,N\in \bN$ with $\max\Set{m,n}< N$. Recall that a random variable $X$ has hypergeometric distribution with parameters $N,n,m$ if $X:=|S\cap [m]|$, where $S$ is a random subset of $[N]$ of size~$n$.
We write $X\sim Bin(n,p)$ if $X$ has binomial distribution with parameters~$n,p$. We will often use the following Chernoff-type bound.

\begin{lemma}[see {\cite[{Corollary~2.3, Remark 2.5, Theorem 2.8 and Theorem 2.10}]{JLR:00}}] \label{lem:chernoff}
Let $X$ be the sum of $n$ independent Bernoulli random variables or let $X$ have a hypergeometric distribution with parameters~$N,n,m$. Then the following hold.
\begin{enumerate}[label={\rm(\roman*)}]
\item For all $t\ge 0$, $\prob{|X - \expn{X}| \geq t} \leq 2\eul^{-2t^2/n}$.\label{chernoff t}
\item For all $0\le\eps \le 3/2$, $\prob{X\neq (1\pm \eps) \expn{X} } \leq 2\eul^{-\eps^2\expn{X}/3}$.\label{chernoff eps}
\end{enumerate}
\end{lemma}

\section{Proof of Theorem~\ref{thm:main}}

\subsection{The final absorbers} \label{absorbers}
In this subsection, we construct the absorbers which will be set aside initially and then used in the final absorption step.
As discussed earlier, our absorbers will consist of the union of `exclusive' absorbers which
can absorb a given graph $L$. Here $L$ is a leftover from a previous
partial decomposition step (i.e.~$L$ plays the role of $H_i$ in the proof overview).

\begin{defin}[Absorber for $L$]
Given a graph~$L$, an \defn{absorber for~$L$} is a graph~$A$ such that $V(L)\In V(A)$ is independent in $A$ and both $A$ and $A\cup L$ have a $K_3$-decomposition.
\end{defin}

Note that the condition that $V(L)$ is independent in $A$ implies that $A$ and $L$ are edge-disjoint.
Observe that if $A$ is an absorber for $L$, then in particular, both $A$ and $A\cup L$ are $K_3$-divisible, and thus $L$ must be $K_3$-divisible. Conversely, we will show that for any $K_3$-divisible graph $L$, there exists an absorber.

To guarantee that we can actually find the constructed absorbers in a given host graph~$G$ of large degree, we will construct absorbers which have low degeneracy.
For a graph~$H$ and a subset $U\In V(H)$, the \defn{degeneracy of $H$ rooted at $U$} is the smallest $d\in \bN\cup\Set{0}$ such that there exists an ordering $v_1,\dots,v_{|H|-|U|}$ of the vertices of $V(H)\sm U$ such that for all $i\in[|H|-|U|]$, $$d_H(v_i,U\cup\set{v_j}{1\le j<i})\le d.$$

The goal of this subsection is to prove the following lemma.

\begin{lemma} \label{lem:absorbers}
Let $L$ be any $K_3$-divisible graph. There exists an absorber $A$ for~$L$ such that the degeneracy of $A$ rooted at $V(L)$ is at most~$4$.
\end{lemma}

We construct absorbers as a concatenation of `transformers', whose purpose is, roughly speaking, to transform a given leftover from a `partial' $K_3$-decomposition into a new leftover. The goal of course is to eventually transform the given leftover $L$ into a new leftover which is $K_3$-decomposable.

\begin{defin}[Transformer]
Given vertex-disjoint graphs~$L,L'$, an \defn{$(L,L')$-transformer} is a graph~$T$ such that $V(L\cup L')\In V(T)$ is independent in $T$ and both $T\cup L$ and $T\cup L'$ have a $K_3$-decomposition.
\end{defin}

Note that the condition that $V(L\cup L')$ is independent in $T$ implies that $T$ is edge-disjoint from both $L$ and~$L'$. Observe also that if $T$ is an $(L,L')$-transformer and $L'$ is $K_3$-decomposable, then $A:=T\cup L'$ is an absorber for~$L$.

Given graphs $H$ and~$H'$, a function $\phi\colon V(H)\to V(H')$ is an \defn{edge-bijective homomorphism from $H$ to $H'$} if $\phi(x)\phi(y)\in E(H')$ for all $xy\in E(H)$, and $e(H)=e(H')=|\set{\phi(x)\phi(y)}{xy\in E(H)}|$.
 We write $H\rightsquigarrow H'$ if such a function $\phi$ exists. More loosely, if $H\rightsquigarrow H'$ then we can merge vertices together in $H$ without creating multi-edges to get a copy of $H'$.

\begin{lemma} \label{lem:transformers}
Let $L$ and $L'$ be vertex-disjoint graphs such that $L\rightsquigarrow L'$ and $2\mid d_L(x)$ for all $x\in V(L)$. There exists an $(L,L')$-transformer~$T$ such that the degeneracy of $T$ rooted at $V(L\cup L')$ is at most~$4$.
\end{lemma}

\proof
Let $\phi\colon L\to L'$ be an edge-bijective homomorphism.
Since $2\mid d_L(x)$ for all $x\in V(L)$, there exists a decomposition $\cC$ of $L$ into cycles. Note that $L'$ decomposes into the graphs $\Set{\phi(C)}_{C\in \cC}$. Suppose that we can find a $(C,\phi(C))$-transformer~$T_C$ for every $C\in \cC$ such that the degeneracy of $T_C$ rooted at $V(C\cup \phi(C))$ is at most~$4$. We may clearly assume that $V(T_C)\cap V(L\cup L')=V(C\cup \phi(C))$ for all $C\in \cC$ and $V(T_C)\cap V(T_{C'}) \In V(L\cup L')$ for all $C,C'\in \cC$. It is then easy to see that $T:=\bigcup_{C\in \cC}T_C$ is an $(L,L')$-transformer such that the degeneracy of $T$ rooted at $V(L\cup L')$ is at most~$4$.

For the remainder of the proof, we may thus assume that $L$ is a cycle $x_1x_2\dots x_sx_1$. We can then construct $T$ as follows. Let $\set{u_i,v_i,w_i}{i\in[s]}$ be a set of $3s$ vertices disjoint from $V(L\cup L')$. The vertex set of $T$ will be $V(L\cup L') \cup \set{u_i,v_i,w_i}{i\in[s]}$. Moreover, we define the following sets of edges (indices modulo~$s$):
\begin{align*}
 E &:= \set{x_iu_i,x_iv_i,x_iw_i,x_iu_{i+1}}{i\in[s]}; 	\\
 E'&:= \set{\phi(x_i)u_i,\phi(x_i)v_i,\phi(x_i)w_i,\phi(x_i)u_{i+1}}{i\in[s]}; 	\\
 \tilde{E} &:= \set{u_iv_i,w_iu_{i+1}}{i\in[s]}; 	\\
 E^\ast &:= \set{v_iw_i}{i\in[s]}.
\end{align*}
Let $E(T):=E\cup E' \cup \tilde{E} \cup E^\ast$ (cf.~Figure~\ref{fig:transformer}).

We claim that $T$ is the desired transformer. Clearly, $V(L\cup L')\In V(T)$ is independent in $T$. Moreover, to see that the degeneracy of $T$ rooted at $V(L\cup L')$ is at most~$4$, order the vertices of $\set{u_i,v_i,w_i}{i\in[s]}$ such that $u_1,\dots,u_s$ come first.
Finally, note that all of $E(L)\cup E \cup E^\ast$,  $E(L')\cup E' \cup E^\ast$, $E\cup \tilde{E}$, $E' \cup \tilde{E}$ have a $K_3$-decomposition. This implies that both $T\cup L$ and $T\cup L'$ have a $K_3$-decomposition.
\endproof

\begin{figure}
    \centering
    \begin{minipage}{0.45\textwidth}
        \centering

\scalebox{0.7}{
\begin{tikzpicture}[
	inner/.style={color=red, thick},
	outer/.style={color=blue, thick},
  pencil/.style={color=gray!20},
  vertex/.style={circle, fill, inner sep=1}]

\coordinate (tr) at (0.6,0.6);
\coordinate (br) at (0.6,-0.6);
\coordinate (tl) at (-0.6,0.6);
\coordinate (bl) at (-0.6,-0.6);

\fill[pencil] (tl) -- (tr) -- (0,3) -- cycle;
\fill[pencil] (bl) -- (br) -- (0,-3) -- cycle;
\fill[pencil] (tl) -- (bl) -- (-3,0) -- cycle;
\fill[pencil] (br) -- (tr) -- (3,0) -- cycle;

\fill[pencil] (tr) -- (1,2) -- (2,1) -- cycle;
\fill[pencil] (br) -- (1,-2) -- (2,-1) -- cycle;
\fill[pencil] (tl) -- (-1,2) -- (-2,1) -- cycle;
\fill[pencil] (bl) -- (-1,-2) -- (-2,-1) -- cycle;

\fill[pencil] (3.8,3.8) -- (3,0) -- (2,1) -- cycle;
\fill[pencil] (3.8,3.8) -- (1,2) -- (0,3) -- cycle;
\fill[pencil] (-3.8,3.8) -- (-3,0) -- (-2,1) -- cycle;
\fill[pencil] (-3.8,3.8) -- (-1,2) -- (0,3) -- cycle;
\fill[pencil] (3.8,-3.8) -- (3,0) -- (2,-1) -- cycle;
\fill[pencil] (3.8,-3.8) -- (1,-2) -- (0,-3) -- cycle;
\fill[pencil] (-3.8,-3.8) -- (-3,0) -- (-2,-1) -- cycle;
\fill[pencil] (-3.8,-3.8) -- (-1,-2) -- (0,-3) -- cycle;

\draw[outer] (3.8,3.8) -- (3.8,-3.8) -- (-3.8,-3.8) -- (-3.8,3.8) -- cycle;
\draw[inner] (tr) -- (br) -- (bl) -- (tl) -- cycle;
\draw (0,3) -- (-3,0) -- (0,-3) -- (3,0) -- cycle;

\foreach \x in {0,...,3}
{
\draw (tr) -- (\x,3-\x);
\draw (tl) -- (-\x,3-\x);
\draw (br) -- (\x,\x-3);
\draw (bl) -- (-\x,\x-3);
\draw (3.8,3.8) -- (\x,3-\x);
\draw (-3.8,3.8) -- (-\x,3-\x);
\draw (3.8,-3.8) -- (\x,\x-3);
\draw (-3.8,-3.8) -- (-\x,\x-3);
}

\node[vertex] at (tr) {};
\node[vertex] at (br) {};
\node[vertex] at (tl) {};
\node[vertex] at (bl) {};
\node at (0.4,0.4) {$x_1$};
\node at (0.4,-0.4) {$x_2$};
\node at (-0.4,-0.4) {$x_3$};
\node at (-0.4,0.4) {$x_4$};

\node[vertex] at (3.8,3.8) {};
\node[vertex] at (-3.8,3.8) {};
\node[vertex] at (3.8,-3.8) {};
\node[vertex] at (-3.8,-3.8) {};
\node at (4.3,3.6) {$\phi(x_1)$};
\node at (4.3,-3.6) {$\phi(x_2)$};
\node at (-4.3,-3.6) {$\phi(x_3)$};
\node at (-4.3,3.6) {$\phi(x_4)$};

\node[vertex] at (3,0) {};
\node at (3.3,0) {$u_2$};
\node[vertex] at (2,1) {};
\node at (1.9,0.7) {$w_1$};
\node[vertex] at (1,2) {};
\node at (0.7,1.9) {$v_1$};
\node[vertex] at (0,3) {};
\node at (0,3.3) {$u_1$};

\node[vertex] at (-3,0) {};
\node at (-3.3,0) {$u_4$};
\node[vertex] at (-2,1) {};
\node at (-1.9,0.7) {$v_4$};
\node[vertex] at (-1,2) {};
\node at (-0.7,1.9) {$w_4$};

\node[vertex] at (2,-1) {};
\node at (1.9,-0.7) {$v_2$};
\node[vertex] at (1,-2) {};
\node at (0.7,-1.9) {$w_2$};
\node[vertex] at (0,-3) {};
\node at (0,-3.4) {$u_3$};

\node[vertex] at (-2,-1) {};
\node at (-1.9,-0.7) {$w_3$};
\node[vertex] at (-1,-2) {};
\node at (-0.7,-1.9) {$v_3$};
\end{tikzpicture}
}

        \caption{A $(C_4, C_4)$-transformer.}
				\label{fig:transformer}
    \end{minipage}\hfill
    \begin{minipage}{0.45\textwidth}
        \centering

\scalebox{0.55}{
\begin{tikzpicture}[
	inner/.style={color=red, thick},
	outer/.style={color=blue, thick},
  pencil/.style={color=gray!20},
  vertex/.style={circle, fill, inner sep=1}]

\coordinate (tr) at (0.6,0.6);
\coordinate (br) at (0.6,-0.6);
\coordinate (tl) at (-0.6,0.6);
\coordinate (bl) at (-0.6,-0.6);

\fill[pencil] (tl) -- (tr) -- (0,3) -- cycle;
\fill[pencil] (bl) -- (br) -- (0,-3) -- cycle;
\fill[pencil] (tl) -- (bl) -- (-3,0) -- cycle;
\fill[pencil] (br) -- (tr) -- (3,0) -- cycle;

\draw[inner] (tr) -- (br) -- (bl) -- (tl) -- cycle;

\draw[outer] (-1,2) -- (0,3) -- (1,2);
\draw[outer] (-1,-2) -- (0,-3) -- (1,-2);
\draw[outer] (2,-1) -- (3,0) -- (2,1);
\draw[outer] (-2,-1) -- (-3,0) -- (-2,1);

\foreach \x in {0,3}
{
\draw[thick] (tr) -- (\x,3-\x);
\draw[thick] (tl) -- (-\x,3-\x);
\draw[thick] (br) -- (\x,\x-3);
\draw[thick] (bl) -- (-\x,\x-3);
}
\foreach \x in {1,2}
{
\draw[outer] (tr) -- (\x,3-\x);
\draw[outer] (tl) -- (-\x,3-\x);
\draw[outer] (br) -- (\x,\x-3);
\draw[outer] (bl) -- (-\x,\x-3);
}
\node[vertex] at (tr) {};
\node[vertex] at (br) {};
\node[vertex] at (tl) {};
\node[vertex] at (bl) {};

\node[vertex] at (3,0) {};
\node[vertex] at (2,1) {};
\node[vertex] at (1,2) {};
\node[vertex] at (0,3) {};

\node[vertex] at (-3,0) {};
\node[vertex] at (-2,1) {};
\node[vertex] at (-1,2) {};

\node[vertex] at (2,-1) {};
\node[vertex] at (1,-2) {};
\node[vertex] at (0,-3) {};

\node[vertex] at (-2,-1) {};
\node[vertex] at (-1,-2) {};

\node at (0,0) {$L$};
\node at (0.0,1.7) {$\nabla L$};
\node at (1.5,2) {$\nabla\nabla L$};

\draw[->, thick,
decorate, decoration={
    zigzag,
    segment length=10,
    amplitude=1.7,post=lineto,
    post length=2pt
}]  (1.7,-1.7) -- (2.2,-2.2);

\begin{scope}[shift={(3.8,-3.8)}]

\coordinate (tr) at (0,0);
\coordinate (br) at (0,0);
\coordinate (tl) at (0,0);
\coordinate (bl) at (0,0);

\fill[pencil] (tl) -- (tr) -- (0,3) -- cycle;
\fill[pencil] (bl) -- (br) -- (0,-3) -- cycle;
\fill[pencil] (tl) -- (bl) -- (-3,0) -- cycle;
\fill[pencil] (br) -- (tr) -- (3,0) -- cycle;

\draw[inner] (tr) -- (br) -- (bl) -- (tl) -- cycle;

\draw[outer] (-1,2) -- (0,3) -- (1,2);
\draw[outer] (-1,-2) -- (0,-3) -- (1,-2);
\draw[outer] (2,-1) -- (3,0) -- (2,1);
\draw[outer] (-2,-1) -- (-3,0) -- (-2,1);

\foreach \x in {1,2}
{
\draw[outer] (tr) -- (\x,3-\x);
\draw[outer] (tl) -- (-\x,3-\x);
\draw[outer] (br) -- (\x,\x-3);
\draw[outer] (bl) -- (-\x,\x-3);
}
\node[vertex] at (tr) {};
\node[vertex] at (br) {};
\node[vertex] at (tl) {};
\node[vertex] at (bl) {};

\node[vertex] at (3,0) {};
\node[vertex] at (2,1) {};
\node[vertex] at (1,2) {};
\node[vertex] at (0,3) {};

\node[vertex] at (-3,0) {};
\node[vertex] at (-2,1) {};
\node[vertex] at (-1,2) {};

\node[vertex] at (2,-1) {};
\node[vertex] at (1,-2) {};
\node[vertex] at (0,-3) {};

\node[vertex] at (-2,-1) {};
\node[vertex] at (-1,-2) {};

\node at (1.2,1.2) {$L_4$};
\end{scope}
\end{tikzpicture}
}
        \caption{The canonical graph $L_4$.}
				\label{fig:canonical}
    \end{minipage}
\end{figure}

We remark that there is a slightly simpler construction of an $(L,L')$-transformer~$T$ if we allowed the degeneracy to be~$6$ (with the effect of obtaining $5/6$ instead of the desired $3/4$ in~Theorem~\ref{thm:main}). Simply add for every edge $e\in E(L)$ a vertex $z_e$ and join $z_e$ to the endvertices of $e$ and $\phi(e)$. Then, for every vertex $x\in V(L)$, add a perfect matching on the vertices $\set{z_{xy}}{y\in N_L(x)}$. This construction turns out to be more easily generalisable to obtain transformers when considering $F$-decompositions rather than $K_3$-decompositions.

As mentioned before, transformers enable us to transform a given leftover from a `partial' $K_3$-decomposition into a new leftover. We now define a `canonical' graph $L_m$ which any $m$-edge $K_3$-divisible graph $L$ can be transformed into (via an intermediary denoted $\nabla\nabla L$). Applying this to both an arbitrary such graph $L$ and a clearly $K_3$-decomposable $m$-edge graph $L'$, and combining the resulting transformers through $L_m$, will give an absorber for $L$.

For $m\in \bN$, let $L_m$ be the graph with
\begin{align*}
	V(L_m) =\Set{v^\ast,v_1,\dots,v_{3m}}, \quad E(L_m) =\bigcup_{i\in[m]}\Set{v^\ast v_{3i-2},v_{3i-2}v_{3i-1},v_{3i-1} v_{3i},v_{3i}v^\ast}.
\end{align*}
That is, $L_m$ is obtained from $C_{4m}$ by taking every fourth vertex on the cycle and identifying all these vertices to a single vertex~$v^\ast$. Note that if $L$ is any graph with $m$ edges, then replacing every edge with a path of length $4$ and identifying the original vertices of $L$ yields a graph isomorphic to~$L_m$ (cf.~Fact~\ref{fact:identifying} and Figure~\ref{fig:canonical}). In this sense, $L_m$ will serve as a `canonical' graph.
Instead of $L_m$, we could have also worked with $C_{4m}$ directly, but this would not generalize to $F$-decompositions for other graphs~$F$, while the current approach does generalize in a natural way.

More formally, given a graph~$L$, we define $\tilde{\nabla} L$ as the graph obtained from $L$ by adding a new set $\set{z_e}{e\in E(L)}$ of vertices disjoint from $V(L)$, and for every edge $e\in E(L)$, we join both endpoints of $e$ to~$z_e$. In other words, we extend every edge of $L$ into a copy of $K_3$. Obviously, $\tilde{\nabla} L$ has a $K_3$-decomposition. Furthermore, let $\nabla L:= \tilde{\nabla} L -L$. Note that $\nabla L$ is obtained from $L$ by replacing every edge with a path of length~$2$, and $\nabla\nabla L$ is obtained from $L$ by replacing every edge with a path of length~$4$.

\begin{fact} \label{fact:identifying}
For every graph $L$ with $m:=e(L)$, we have $\nabla\nabla L \rightsquigarrow L_{m}$.
\end{fact}

\proof
Define $\phi\colon \nabla\nabla L \to  L_{m}$ such that $\phi(x):=v^\ast$ for all $x\in V(L)$, and map $V(\nabla\nabla L)\sm V(L)$ bijectively to $V(L_m)\sm \Set{v^\ast}$ in the obvious way.
\endproof

We now construct the absorber $A$ as the union of several intermediate graphs and transformers.
\lateproof{Lemma~\ref{lem:absorbers}}
Let $m:=e(L)$. Let $L'$ be the vertex-disjoint union of $m/3$~triangles. Moreover, let $L_m$ be the `canonical' graph defined above, and also define $\nabla L$, $\nabla\nabla L$, $\nabla L'$, $\nabla\nabla L'$ as above. We assume that $\nabla\nabla L$, $\nabla\nabla L'$ and $L_m$ are pairwise vertex-disjoint.

By Fact~\ref{fact:identifying}, we have $\nabla\nabla L \rightsquigarrow L_m$ and $\nabla\nabla L' \rightsquigarrow L_m$.
 Thus, by Lemma~\ref{lem:transformers}, there exists an $(\nabla\nabla L,L_m)$-transformer~$T$ such that the degeneracy of $T$ rooted at $V(\nabla\nabla L\cup L_m)$ is at most~$4$, and there exists an $(\nabla\nabla L',L_m)$-transformer~$T'$ such that the degeneracy of $T'$ rooted at $V(\nabla\nabla L'\cup L_m)$ is at most~$4$. We may assume that $T$ and $T'$ consist of new vertices except for the unavoidable intersections, that is, $V(T)\cap V(T')=V(L_m)$, $V(T)\cap V(\nabla\nabla L')=\emptyset$, $V(T')\cap V(\nabla\nabla L)=\emptyset$.

We now define the graph
\begin{align*}
	A:= \nabla L \cup \nabla \nabla L \cup T \cup L_m \cup T' \cup \nabla \nabla L' \cup \nabla L' \cup L'
\end{align*}
and claim that $A$ is the desired absorber. Clearly, $V(L)$ is independent in~$A$, and it is easy to see that the degeneracy of $A$ rooted at $V(L)$ is at most $4$.\COMMENT{Order the vertices of $V(A)\sm V(L)$ as follows: take the distinguished vertex $v^\ast$ of $V(L_m)$ first, followed by $V(L_m)\sm \Set{v^\ast}$, then continue with $V(L')$ in any order, $V(\nabla\nabla L)\sm V(L)$, $V(\nabla\nabla L')\sm V(L')$ in any order; ultimately, order the vertices of $V(T)\sm V(\nabla\nabla L \cup L_m)$ and of $V(T') \sm V(\nabla\nabla L'\cup L_m)$ in such a way that every vertex is connected to at most~$4$ previous vertices.}

Finally, note that $\nabla L \cup \nabla \nabla L=\tilde{\nabla}(\nabla L)$, $T \cup L_m$, $T' \cup \nabla\nabla L'$ and $\nabla L' \cup L'=\tilde{\nabla}L'$ are pairwise edge-disjoint and are all $K_3$-decomposable. Thus, $A$ has a $K_3$-decomposition. Similarly,
\begin{align*}
	A\cup L &= (L\cup \nabla L) \cup (\nabla\nabla L \cup T) \cup (L_m \cup T') \cup (\nabla \nabla L' \cup \nabla L') \cup L'
\end{align*}\COMMENT{$= \tilde{\nabla}L \cup (\nabla\nabla L \cup T) \cup (L_m \cup T') \cup \tilde{\nabla} (\nabla L') \cup L'$}
 has a $K_3$-decomposition.
\endproof

\subsection{Vortices}
Our proof method involves an iterative absorption process, where in each iteration, the number of leftover configurations is drastically reduced. We ensure this by forcing leftover edges into smaller and smaller subsets of the vertex set. The underlying structure is a `vortex', which is defined as follows.

\begin{defin}[Vortex]\label{def:vortex}
Let $G$ be a graph on $n$ vertices. A \defn{$(\delta,\eps,m)$-vortex in $G$} is a sequence \mbox{$U_0 \supseteq U_1 \supseteq \dots \supseteq U_\ell$} such that
\begin{enumerate}[label=\rm{(V\arabic*)}]
\item $U_0=V(G)$; \label{vortex 1}
\item $|U_i|=\lfloor \eps|U_{i-1}| \rfloor$ for all $i \in [\ell]$; \label{vortex 2}
\item $|U_\ell|=m$; \label{vortex 3}
\item $d_G(x,U_i) \ge \delta |U_i|$ for all $i \in [\ell]$ and $x \in U_{i-1}$. \label{vortex 4}
\end{enumerate}
\end{defin}

We now show that every large graph of high minimum degree contains a vortex such that the final set $U_\ell$ has constant size.
This can easily be achieved by successively taking random subsets.

\begin{lemma} \label{lem:get vortex}
Let $\delta\in[0,1]$ and $1/m'\ll \eps<1$. Suppose that $G$ is a graph on $n \ge m'$ vertices with $\delta(G) \ge \delta n$. Then $G$ has a $(\delta-\eps,\eps,m)$-vortex for some $\lfloor \eps m' \rfloor \le m \le m'$.
\end{lemma}

\proof
Recursively, define  $n_0:=n$ and $n_i:=\lfloor \eps n_{i-1} \rfloor$. Observe that $\eps^i n \ge n_i \ge \eps^i n- 1/(1-\eps)$.\COMMENT{$n_i \ge \eps^i n -1-\eps-\ldots -\eps^{i-1}\geq \eps^i n- 1/(1-\eps)$} Let $\ell:=1+\max\set{i\ge 0}{n_i\ge m'}$ and let $m:=n_\ell$. Note that $\lfloor \eps m' \rfloor \le m \le m'$. Further, for $i\in[\ell]$, let
\begin{align}\label{newone}\eps_i:=n^{-1/3}\sum_{j=1}^i \eps^{-(j-1)/3}=n^{-1/3}\frac{\eps^{-i/3}-1}{\eps^{-1/3}-1}  \le \frac{(\eps^{i-1}n)^{-1/3}}{1-\eps^{1/3}} \le \frac{m'^{-1/3}}{1-\eps^{1/3}}  \le \frac{\eps}{3},
\end{align} where we have used $\eps^{i-1}n \ge \eps^{\ell-1}n\ge n_{\ell-1}\ge m'$,
 and let $\eps_0:=0$.

Now, suppose that for some $i\in[\ell]$, we have already found a $(\delta-3\eps_{i-1},\eps,n_{i-1})$-vortex $U_0,\dots,U_{i-1}$ in $G$ (which is true for $i=1$). In particular, $\delta(G[U_{i-1}]) \ge (\delta-3\eps_{i-1})n_{i-1}$.
Choose $U_i$ to be a random subset of~$U_{i-1}$ of size~$n_i$. Lemma~\ref{lem:chernoff} implies that with positive probability, $d_G(x,U_i) \ge (\delta-3\eps_{i-1}-2n_{i-1}^{-1/3})n_i$ for all $x\in U_{i-1}$.\COMMENT{For every $x\in U_{i-1}$, $$\prob{d_G(x,U_i)<(\delta -3\eps_{i-1}) \lfloor \eps n_{i-1} \rfloor -\eps n_{i-1}^{2/3}} \le 2 \eul^{-2\eps^2 n_{i-1}^{4/3}/\lfloor \eps n_{i-1}\rfloor}.$$} Fix such a choice of~$U_i$.
Then $U_0,\dots,U_i$ form a $(\delta-3\eps_i,\eps,n_i)$-vortex in~$G$.
Repeating this for all $i\in[\ell]$, we finally obtain a $(\delta-3\eps_\ell,\eps,m)$-vortex $U_0,\dots,U_\ell$ in $G$. As, by~\eqref{newone}, $\eps_\ell\le \eps/3$, the lemma follows.
\endproof

\subsection{Cover down lemma}
The engine behind the iterative absorption procedure is the following `Cover down lemma'. Recall that the definition of $\delta^{0+}$ allows us to find approximate $K_3$-decompositions such that the maximum degree of the leftover is very small. The strength of the Cover down lemma is that it also forces all leftover edges to lie inside a prescribed small vertex set~$U$ (which one might think of as the next vortex set~$U_i$ in the sequence).

\begin{lemma}[Cover down lemma] \label{lem:cover-down}
Suppose $1/n \ll \eps$ and let $\delta:=\max\Set{2/3,\delta^{0+}}$. Let $G$ be a graph on $n$ vertices and $U\In V(G)$ with $|U|=\lfloor \eps n \rfloor$.
Suppose that $\delta(G) \ge (\delta+3\eps)n$ and $d_G(x,U)\ge (\delta + 2\eps)|U|$ for all $x\in V(G)$. Also assume that $d_G(x)$ is even for all $x\in V(G)\sm U$. Then there exists a $K_3$-decomposable subgraph $H\In G$ such that $G-G[U] \In H$ and $\Delta(H[U])\le \eps^{10} n$.
\end{lemma}

In order to prove Lemma~\ref{lem:cover-down}, we will first randomly select a sparse reserve graph before using the definition of $\delta^{0+}$ to find an approximate $K_3$-decomposition of $G-G[U]$. We then cover all leftover edges which are not inside~$U$ in two stages. First, we find for every leftover edge $e=xy$ with $x,y\in V(G)\sm U$ a vertex $u_e\in U$ such that $u_e,x,y$ form a triangle in~$G$. In a second step, we cover the remaining `crossing' leftover edges. Suppose $x\in V(G)\sm U$ and let $U_x$ be the set of all remaining neighbours of~$x$.
What we would now like to find is a perfect matching $M_x$ of the `link graph' $G[U_x]$ of~$x$. Clearly, the edges of $M_x$ together with all the edges from $U_x$ to $x$ would then form edge-disjoint triangles covering all remaining edges at~$x$. Lemma~\ref{lem:pseudorandom remainder} will allow us to achieve the latter simultaneously for all $x\in V(G)\sm U$. For this, we need the following simple result.

\begin{fact}[cf.~{\cite[Lemma 8]{raman:90}}] \label{fact:Jain}
Let $X_1, \ldots, X_n$ be Bernoulli random variables such that for all $i \in [n]$, we have $\prob{X_i = 1 \mid X_1, \ldots, X_{i-1}} \leq p$.
Let $B \sim Bin(n,p)$ and $X:=\sum_{i=1}^n X_i$. Then $\prob{X \geq a} \leq \prob{B \geq a}$ for all~$a\ge 0$.
\end{fact}

\begin{lemma} \label{lem:pseudorandom remainder}
Let $1/n \ll \rho$ and $N\in \bN$. Let $G$ be a graph on $n$ vertices and suppose there are sets $U_1,\dots,U_N \In V(G)$ with the following properties:
\begin{enumerate}[label=(\roman*)]
\item $2 \mid |U_i|$ and $\delta(G[U_i]) \ge (1/2+4\rho^{1/6})|U_i|$ for all $i \in [N]$; \label{eqn:coverUW:degrees and div}
\item $|U_i| \ge \rho^{4/3} n$ for all $i\in[N]$; \label{eqn:coverUW:large enough}
\item $|U_i \cap U_j| \le \rho^2 n$ for all $1\le i < j\le N$; \label{eqn:coverUW:intersection}
\item every vertex $u\in V(G)$ is contained in at most $\rho n$ of the sets $U_i$. \label{eqn:coverUW:spread}
\end{enumerate}
Then for every $i\in [N]$, there exists a perfect matching $M_i$ of $G[U_i]$, such that all the matchings $\Set{M_i}_{i\in[N]}$ are pairwise edge-disjoint.
\end{lemma}

It follows directly from \ref{eqn:coverUW:degrees and div} that for every $i\in [N]$, there exists a perfect matching $M_i$ of $G[U_i]$. The difficulty here lies in finding edge-disjoint ones. For this, we use a randomised algorithm.

\proof
Let $t: = \lceil 2 \rho^{3/2} n \rceil $ and define $G_i:=G[U_i]$ for all $i\in [N]$.
Suppose that we have already found $M_1, \dots, M_{i-1}$ for some $i \in [N]$. We now define $M_i$ as follows.
Let $H_{i-1} : = \bigcup_{j=1}^{i-1}  M_{j}$ and let $G_i' : = (G - H_{i-1})[U_i]$.
If $\Delta ( H_{i-1}[U_i] ) \le \rho^{3/2} n$, then
\begin{align*}
\delta(G'_i) &\ge \delta(G[U_i])-\Delta(H_{i-1}[U_i]) \geq (1/2+4\rho^{1/6})|U_i| -  \rho^{3/2} n \ge |U_i|/2 + t
\end{align*}
by \ref{eqn:coverUW:degrees and div} and \ref{eqn:coverUW:large enough}.\COMMENT{$4\rho^{1/6}|U_i|\ge 4\rho^{1/6}\rho^{4/3} n \ge t+\rho^{3/2} n$} Thus, we can successively find $t$ edge-disjoint perfect matchings $A_1,\dots,A_t$ of~$G_i'$ (which are all suitable candidates for~$M_i$). Otherwise, if $ \Delta ( H_{i-1}[U_i] ) >  \rho^{3/2} n$, then let $A_1, \ldots, A_t$ be empty graphs on~$U_i$.

In either case, we have found edge-disjoint subgraphs $A_1, \dots, A_t$ of $G_i'$.
Pick $s \in [t]$ uniformly at random and set $M_i : = A_s$.
The lemma follows if the following holds with positive probability:
\begin{align}
\Delta ( H_{i-1}[U_i] ) \le   \rho^{3/2} n \text{ for all $i\in[N]$.}	\label{eqn:Krmkey}
\end{align}
For $i \in [N]$ and $u \in U_i$, let $J^{i,u}$ be the set of indices $j\in [i-1]$ such that $u \in U_{j}$, and for $j\in J^{i,u}$, let $Y^{i,u}_j$ be the indicator variable of the event that $uu'\in E(M_j)$ for some $u'\in U_i$.
Observe that $$d_{H_{i-1}[U_i]}(u)=\sum_{j\in J^{i,u}}Y^{i,u}_j.$$

Now, fix $i \in [N]$ and $u \in U_i$.
Crucially, for any $j\in J^{i,u}$, by~\ref{eqn:coverUW:intersection}, at most $\rho^2 n$ of the subgraphs $A_1\dots,A_t$ that we picked in $G'_j$ contain an edge incident to~$u$ in~$G_i$ (regardless of the previous choices).
Let $j_1, \dots, j_{|J^{i,u}|}$ be the enumeration of $J^{i,u} $ in increasing order. By the above, for all $\ell \in [|J^{i,u}|] $, we have
\begin{align*}
	\prob{ Y^{i,u}_{j_\ell} = 1 \mid Y^{i,u}_{j_1}, \dots, Y^{i,u}_{j_{\ell-1}}} \le
	\frac{\rho^2 n}{ t } \le \frac {\rho^{1/2}} {2}.
\end{align*}
Let $B \sim Bin( |J^{i,u}|  , \rho^{1/2}/2 )$. Since $ |J^{i,u}|  \le \rho n$ by~\ref{eqn:coverUW:spread}, we have $\expn{B}\le \rho^{3/2}n/2$.
Using Fact~\ref{fact:Jain} and Lemma~\ref{lem:chernoff}, we infer that
\begin{align*}
\prob{ \sum_{j \in J^{i,u}} Y^{i,u}_{j}  >  \rho^{3/2} n }
\le \prob{ B >  \rho^{3/2} n }
 \le \prob{ B > \expn{B} +  \rho^{3/2} n/2 }
\le 2 e^{-  \rho^{2} n/2}.
\end{align*}

Finally, since by~\ref{eqn:coverUW:spread} there are at most $\rho n^2$ pairs $(i, u)$ with $u\in U_i$, a union bound implies that~\eqref{eqn:Krmkey} holds with positive probability.
\endproof

\lateproof{Lemma~\ref{lem:cover-down}}
Choose new constants $\gamma,\rho>0$ such that $1/n\ll \gamma \ll \rho \ll \eps$. Let $W:=V(G)\sm U$ and let $w_1,\dots,w_N$ be an enumeration of~$W$.

We first observe that since $d_G(x,U)\ge (2/3 + 2\eps)|U|$ for all $x\in V(G)$, it follows that for all $x,y\in V(G)$, we have
\begin{align}
	|N_G(x)\cap N_G(y)\cap U| \ge (1/2+3\eps)|N_{G}(y)\cap U| \ge |U|/3 .  \label{link intersection}
\end{align}\COMMENT{Let $X:=N_G(x)\cap U$ and $Y:=N_G(y)\cap U$. Then $|X|+|Y|\le |U|+|X\cap Y|$ and hence $|X\cap Y| \ge |X|+|Y|-|U| \ge |Y| + (2/3+2\eps)|U|-|U| \ge |Y|-(1/3-2\eps)|U|$. With $|U|\le 3|Y|/2$, it follows that $|X\cap Y| \ge (1/2+3\eps) |Y|$ }

Before obtaining an approximate decomposition, we set aside a sparse graph~$R$ which will act as a `partial absorber'.
For this, let $U_1,\dots,U_N$ be sets with the following properties:
\begin{enumerate}[label=\rm{(\alph*)}]
\item $U_i\In N_G(w_i)\cap U$ for all $i\in[N]$;      \label{eqn:reserve1:neighbourhood}
\item $|U_i|=(1\pm \rho)\rho |N_G(w_i)\cap U|$ for all $i\in[N]$;  \label{eqn:reserve1:size}  
\item $\rho^2 |U|/4\le |U_i \cap U_j| \le 2\rho^2 |U|$ for all $1\le i < j\le N$;  \label{eqn:reserve1:intersection}
\item $|N_G(u)\cap U_i| \ge (1-\rho)\rho (1/2+3\eps)|N_G(w_i)\cap U|$ for all $u\in U$ and $i\in[N]$;\label{eqn:reserve1:degrees}
\item each $u\in U$ is contained in at most $2\rho n$ of the $U_i$'s.           \label{eqn:reserve1:spread}
\end{enumerate}
That such subsets exist can be seen via a probabilistic argument. Indeed, for every pair $(u,i)$ with $i\in[N]$ and $u\in N_G(w_i)\cap U$, include $u$ in $U_i$ with probability~$\rho$ (independently of all other pairs). Applying Lemma~\ref{lem:chernoff} (and using~\eqref{link intersection}) shows that the random sets $U_1,\dots,U_N$ satisfy the desired properties with positive probability.\COMMENT{\ref{eqn:reserve1:neighbourhood} clearly holds. \ref{eqn:reserve1:size} $\expn{|U_i|}=\rho d_G(w_i,U)$, where $|U|\ge d_G(w_i,U)\ge 2|U|/3$. \ref{eqn:reserve1:intersection} $\expn{|U_i\cap U_j|}=\rho^2 |N_G(w_i)\cap N_G(w_j)\cap U|$, where $|U|\ge |N_G(w_i)\cap N_G(w_j)\cap U| \ge |U|/3$. \ref{eqn:reserve1:spread} expected number is $\rho N\le \rho n$. \ref{eqn:reserve1:degrees} $d_G(u,N_G(w_i)\cap U)\ge (1/2+3\eps)|N_{G}(w_i)\cap U|$.}
Note that it follows from~\ref{eqn:reserve1:size} and~\ref{eqn:reserve1:degrees} that for all $u\in U$ and $i\in [N]$, we have
\begin{align}
	d_G(u,U_i)\ge (1/2+2\eps)|U_i|.\label{eqn:reserve1:degrees new}
\end{align}
Let $R$ be the subgraph of $G[U,W]$ consisting of all edges $uw_i$ where $i\in[N]$ and $u\in U_i$. Note that $\Delta(R)\le 2\rho n$ by~\ref{eqn:reserve1:size} and \ref{eqn:reserve1:spread}.

Let $G':=G-G[U]-R$. Clearly, we have $\delta(G')\ge (\delta^{0+}+\eps)n$. By definition of $\delta^{0+}$, there exists a subgraph $L$ of $G'$ such that $\Delta(L)\le \gamma n$ and $G'-L$ is $K_3$-decomposable.

Next, for every edge $e=w_iw_j \in E(L[W])$, we choose a vertex $u_e \in U_i\cap U_j$ in such a way that $u_e\neq u_{e'}$ whenever $e\cap e'\neq \emptyset$.

This can be done greedily. Indeed, whenever we want to choose~$u_{w_iw_j}$, there are at least $\rho^2 |U|/4$ vertices $u\in U_i\cap U_j$ by~\ref{eqn:reserve1:intersection}. Moreover, at most $2\Delta(L)$ of these vertices $u$ are blocked by some edge $e'$ which has $w_i$ or $w_j$ as an endpoint and $u_{e'}=u$ has been previously chosen. Since $2\Delta(L) < \rho^2 |U|/4$, we can always choose a suitable vertex~$u_{w_iw_j}$.

Let $\hat{H}$ be the graph consisting of all the edges $u_{w_iw_j}w_i,u_{w_iw_j}w_j$ with $w_iw_j \in E(L[W])$. By~\ref{eqn:reserve1:neighbourhood}, we can see that $\hat{H}$ is a subgraph of~$R$. Moreover, $\hat{H}\cup L[W]$ clearly has a $K_3$-decomposition.
Let $$R':= (R- \hat{H}) \cup L[U,W].$$ By the above, $G-G[U]-R'=(G'-L)\cup (\hat{H}\cup L[W])$ has a $K_3$-decomposition. It remains to cover all the edges of $R'$ using only a few edges of~$G[U]$.
For every $i\in[N]$, let $$U_i':= N_{R'}(w_i) = (U_i\sm N_{\hat{H}}(w_i)) \cup N_{L[U,W]}(w_i).$$ Since $G[U]\cup R'$ is obtained from $G$ by removing edge-disjoint triangles, we have that $|U_i'|$ is even for all $i\in[N]$.
Moreover, since $d_{\hat{H}}(w_i)=d_{L[W]}(w_i) \le \gamma n$ and $d_{L[U,W]}(w_i)\le \gamma n$ for all $i\in[N]$, we have $|U_i'|\ge \rho|U|/2$ for all $i\in[N]$ and $|U_i'\cap U_j'|\le 3\rho^2 |U|$ for all $1\le i<j \le N$. Since $d_{L[U,W]}(u)\le \gamma n$ for all $u\in U$, it also follows that every vertex $u\in U$ is contained in at most $3\rho n \le 4\rho \eps^{-1}|U|$ of the $U_i'$'s. Moreover, we can deduce from~\eqref{eqn:reserve1:degrees new} that $\delta(G[U_i']) \ge (1/2+\eps)|U_i'|$ for all $i \in [N]$.
Thus, by Lemma~\ref{lem:pseudorandom remainder} (with $G[U]$, $|U|,4\rho/\eps$ playing the roles of $G,n,\rho$), for every $i\in[N]$, there exists a perfect matching of $G[U_i']$, such that all the matchings $\Set{M_i}_{i\in[N]}$ are pairwise edge-disjoint.
Then
\begin{align*}
	\bigcup_{i\in[N]}M_i \cup R' = \bigcup_{i\in[N]} (M_i \cup \set{w_i u}{u\in N_{R'}(w_i)})
\end{align*}
is $K_3$-decomposable. Thus $H:=(G-G[U])\cup \bigcup_{i\in[N]}M_i$ is $K_3$-decomposable. Moreover, $\Delta(H[U])=\Delta(\bigcup_{i\in[N]}M_i) \le 3\rho n\le \eps^{10} n$.
\endproof

\subsection{Proof of Theorem~\ref{thm:main}}

We can now combine Lemmas~\ref{lem:absorbers},~\ref{lem:get vortex} and~\ref{lem:cover-down} to prove Theorem~\ref{thm:main}.

\lateproof{Theorem~\ref{thm:main}}
For convenience, we will prove this with $8\eps$ in place of $\eps$. That is, we assume that $\delta(G)\ge (\delta +8\eps)n$.
Choose new constants $m',M\in \bN$ such that $1/n_0 \ll  1/M \ll 1/m' \ll \eps \ll 1$. Let $\delta:= \max\Set{\delta^{0+},3/4}$. By Lemma~\ref{lem:absorbers}, and as $1/M \ll 1/m'$, for every $K_3$-divisible graph $L$ with $|L|\le m'$, there exists an absorber $A_L$ for~$L$ such that $|A_L|\le M$ and the degeneracy of $A_L$ rooted at $V(L)$ is at most~$4$.

Let $G$ be a $K_3$-divisible graph on $n\ge n_0$ vertices with $\delta(G)\ge (\delta +8\eps)n$. Our aim is to show that $G$ has a $K_3$-decomposition. We achieve this in four steps.

{\bf Step 1.} First, we apply Lemma~\ref{lem:get vortex} to obtain a $(\delta+7\eps,\eps,m)$-vortex $U_0,U_1,\dots,U_\ell$ in $G$ for some $\lfloor \eps m' \rfloor \le m \le m'$.

{\bf Step 2.} Next, we find `exclusive' absorbers for the possible leftover graphs on~$U_\ell$.
To this end, let $\cL$ be the collection of all spanning $K_3$-divisible subgraphs of~$G[U_\ell]$. Obviously, $|\cL|\le 2^{\binom{m}{2}}$. It is thus easy to find edge-disjoint subgraphs $\Set{\tilde{A}_L}_{L\in \cL}$ of~$G$ such that for all $L\in\cL$, we have that $\tilde{A}_L$ is an absorber for~$L$, $|\tilde{A}_L|\le M$, and $\tilde{A}_L[U_1]$ is empty. Indeed, we can find these graphs in turn. Suppose we want to find $\tilde{A}_L$. Consider the graph $\tilde{G}$ obtained from $G-G[U_1]$ by deleting the edges of previously chosen absorbers. Note that $\delta(\tilde{G})\ge (3/4+\eps)n$. Thus, any four vertices in $\tilde{G}$ have at least $4\eps n$ common neighbours. Since there is an ordering of the vertices of $V(A_L)\sm V(L)$ such that every vertex is joined to at most $4$ preceding vertices in~$A_L$,\COMMENT{by the above, where Lemma~\ref{lem:absorbers} is already applied} we can embed the vertices of $V(A_L)\sm V(L)$ one after the other into $\tilde{G}$ to obtain~$\tilde{A}_L$.

Let $A^\ast:=\bigcup_{L\in \cL} \tilde{A}_L$. Observe that $A^\ast$ has the following crucial property:
\begin{align}
	\mbox{given any $K_3$-divisible subgraph $L^\ast$ of $G[U_\ell]$, $A^\ast\cup L^\ast$ has a $K_3$-decomposition.} \label{absorber all:property}
\end{align}
Let $G':=G-A^\ast$.
Observe that since $\Delta(A^\ast)\le M |\cL|$ and $A^\ast[U_1]$ is empty, we have that $U_0,U_1,\dots,U_\ell$ is a $(\delta+6\eps,\eps,m)$-vortex in~$G'$ and $\delta(G')\ge (\delta+7\eps)n$. Moreover, since $A^\ast$ is the edge-disjoint union of absorbers, it must be $K_3$-divisible, and thus $G'$ is also $K_3$-divisible.

{\bf Step 3.} We now iteratively apply the Cover down lemma (Lemma~\ref{lem:cover-down}) to cover all the edges of $G'$ except possibly some inside~$U_\ell$.
More precisely, we show inductively that for all $i\in[\ell]\cup \Set{0}$, there exists a subgraph $G_i \In G'[U_i]$
such that $G'-G_i$ has a $K_3$-decomposition, and such that the following hold (where $U_{\ell+1}:=\emptyset$):
\begin{align}
	\delta(G_i) &\ge (\delta+4\eps)|U_i|;  \label{cover down induction 1}  \\
	d_{G_i}(x,U_{i+1}) &\ge (\delta + 5\eps)|U_{i+1}| \mbox{ for all }x\in U_{i}; \label{cover down induction 2}  \\
	G_i[U_{i+1}] &=G'[U_{i+1}]. \label{cover down induction 3}
\end{align}
Clearly, this holds for $i=0$ with $G_0:=G'$. Now, suppose that for some $i\in[\ell-1]\cup \Set{0}$, we have found $G_i$ satisfying the above. Note that $G_i$ is $K_3$-divisible. Define $G_i':=G_i-G_i[U_{i+2}]$. We still have that $d_{G_i'}(x)$ is even for all $x\in U_i\sm U_{i+1}$.
Thus, by Lemma~\ref{lem:cover-down} (with $G_i'$, $U_{i+1}$, $\eps$ playing the roles of $G,U,\eps$),\COMMENT{$\delta(G_i')\ge (\delta+3\eps)|U_i|$ and $d_{G_i'}(x,U_{i+1})\ge (\delta+4\eps)|U_{i+1}|$ for all $x\in U_i$. $\eps$ has to play the role of $\eps$ because it is also the ratio of $|U_{i+1}|/|U_i|$.} there exists a $K_3$-decomposable subgraph $H\In G_i'$ such that
$G_i'-G_i'[U_{i+1}] \In H$ and $	\Delta(H[U_{i+1}])\le \eps^8 |U_{i+1}|$.
Let $G_{i+1}:=(G_i-H)[U_{i+1}].$

Note that since $G_i-G_i[U_{i+1}]=G_i'-G_i'[U_{i+1}] \In H$,
we can deduce that $G'-G_{i+1}=(G'-G_i) \cup H$\COMMENT{Clearly, $H\In G'$. Suppose $e\in G'-G_{i+1}$. Then we must have $e\notin G_i[U_{i+1}]$ or $e\in H$. In the latter case, we're done. In the first case, if $e\in G_i$, we can deduce $e\in H$ as well using the previous identity, and if $e\notin G_i$, we're done as well. For the converse, suppose $e\in G'-G_i$ or $e\in H$. This immediately gives $e\in G'-G_{i+1}$.} is $K_3$-decomposable, as desired.
Moreover, since $H[U_{i+2}]$ is empty by definition of~$G_i'$, we clearly have $G_{i+1}[U_{i+2}] =G'[U_{i+2}]$. Observe that
\begin{align*}
	\delta(G_{i+1}) \overset{\eqref{cover down induction 2}}{\ge} (\delta+5\eps)|U_{i+1}| - \Delta(H[U_{i+1}])  \ge   (\delta+4\eps)|U_{i+1}|
\end{align*}\COMMENT{for first inequality, can also use \eqref{cover down induction 3},\ref{vortex 4}.}
and for every $x\in U_{i+1}$, we have
\begin{align*}
  d_{G_{i+1}}(x,U_{i+2})  \overset{\eqref{cover down induction 3},\ref{vortex 4}}{\ge} (\delta+6\eps)|U_{i+2}| - \Delta(H[U_{i+1}])  \ge  (\delta + 5\eps)|U_{i+2}|.
\end{align*}
Thus, \eqref{cover down induction 1}--\eqref{cover down induction 3} hold with $i$ replaced by~$i+1$.
By induction, there exists a subgraph $G_\ell \In G'[U_\ell]$ such that $G'-G_\ell$ has a $K_3$-decomposition.

{\bf Step 4.} Finally, since $G_\ell$ is $K_3$-divisible, $A^\ast\cup G_\ell$ has a $K_3$-decomposition by~\eqref{absorber all:property}. Altogether, $G=(G'-G_\ell) \cup (A^\ast\cup G_\ell)$ has a $K_3$-decomposition, as desired.
\endproof

\section{Proof of Theorem~\ref{thm:qr}} \label{sec:qr}

The proof of Theorem~\ref{thm:main} can be easily adapted to prove Theorem~\ref{thm:qr}, modulo a `boosting step' which we discuss below. In particular, the construction of the absorbers does not need any changes. Obviously, the definition of a vortex has to be adapted to the quasirandom setting. More precisely, instead of condition \ref{vortex 4} in Definition~\ref{def:vortex}, one now requires that for every set $A \subseteq U_{i-1}$ with $|A| \le 4$ the common neighbourhood in $U_i$ of the vertices in $A$ has size $(1\pm \xi)p^{|A|}|U_i|$, where $\xi$ is a small error parameter. This definition ensures that the typicality condition is preserved throughout the iterative absorption procedure. Moreover, it allows to perform the Cover down step for each~$i\in[\ell]$. In this step, we now use that for every vertex $x\in U_{i-1}$, the link graph $G[N_G(x)\cap U_i]$ will be $(\sqrt{\xi},3,p)$-typical, replacing the minimum degree condition in Lemma~\ref{lem:pseudorandom remainder}\ref{eqn:coverUW:degrees and div}. It is well known that typicality ensures the existence of perfect matchings as needed, e.g.~one could split the graph randomly into two equal-sized parts, use the fact (see~\cite{DLR:95}) that typicality implies super-regularity, and check Hall's condition.

Theorem~\ref{thm:qr} has the advantage that it does not rely on the (unknown) parameter~$\delta^{0+}$. Indeed, in this setting, the required approximate decompositions can be obtained using the nibble method, introduced by R\"odl~\cite{rodl:85}, which has since had an enormous influence on combinatorics.
A collection $\cT$ of triangles in a graph $G$ is \defn{$(\xi,p)$-regular} if every edge of $G$ is contained in $(1\pm \xi)p^2 n$ triangles.
The following follows for instance from a result in~\cite{AY:05}.

\begin{cor} \label{cor:qr approx}
Let $1/n\ll \xi \ll \gamma,p$. Let $G$ be a graph on $n$ vertices which contains a $(\xi,p)$-regular collection of triangles. Then $G$ contains a $K_3$-decomposable subgraph $H$ such that $\Delta(G-H)\le \gamma n$.
\end{cor}

Clearly, the collection of all triangles of a $(\xi,2,p)$-typical graph is $(\xi,p)$-regular, thus we immediately obtain approximate $K_3$-decompositions in such a graph.
However, Corollary~\ref{cor:qr approx} itself is not sufficient to replace $\delta^{0+}$ in the proof of Theorem~\ref{thm:main}. This is because of the requirement $\xi \ll \gamma$. Though this is a reasonable assumption in Corollary~\ref{cor:qr approx}, it would make it impossible to control the error parameter $\xi$ during the iterative absorption procedure, since the parameter $\gamma$, which controls the maximum degree of the leftover of the approximate decomposition, feeds (via the Cover down step) into the typicality parameter for the subsequent iteration step.
To overcome this, we need an intermediate step which `boosts' the regularity parameter. The aim is to guarantee approximate decompositions with leftover maximum degree $\gamma n$ in a $(\xi,4,p)$-typical graph even if $\gamma \ll \xi$. To achieve this, instead of applying Corollary~\ref{cor:qr approx} to the collection of all triangles of~$G$, we find a suitable subcollection which is $(\xi',p)$-regular, where $\xi'\ll \gamma$, to which we can then apply Corollary~\ref{cor:qr approx}. The idea is to choose such a collection randomly, according to a suitable probability distribution. Such a probability distribution can be thought of as a `fractional triangle-equicovering', which is similar to a fractional triangle decomposition. Indeed, our basic tool to find such a fractional triangle-equicovering is a so-called \defn{$K_5$-shifter}, which was introduced in the context of fractional decompositions in~\cite{BKLMO:17}. A $K_5$-shifter is a `local' function which allows to adjust the total weight of one edge without affecting the total weight of any other edge. A similar idea of `regularity boosting' has also been successfully applied e.g.~in~\cite{AKS:97}.

\begin{lemma} \label{lem:boosting}
Let $1/n\ll \xi ,p$ and $\xi\le p^7/20$. Any $(\xi,4,p)$-typical graph on $n$ vertices contains an $(n^{-1/3},p/2)$-regular collection of triangles.
\end{lemma}

\proof
Let $G$ be a $(\xi,4,p)$-typical graph on $n$ vertices, and let $\cT^{(3)}$ be the collection of all triangles in~$G$. Moreover, let $\cT^{(5)}$ be the collection of all $K_5$'s in~$G$. For every $e\in E(G)$ and $i\in\Set{3,5}$, let $\cT^{(i)}(e)$ denote the set of all elements of $\cT^{(i)}$ which contain~$e$ as an edge.

Assume, for the moment, that $\psi\colon \cT^{(3)}\to [0,1]$ is a function such that for every edge $e\in E(G)$, we have $$\sum_{T\in \cT^{(3)}(e)}\psi(T)=\frac{1}{4}p^2 n.$$ We can then choose a random subcollection $\cT'\In \cT^{(3)}$ by including every $T\in \cT^{(3)}$ with probability $\psi(T)$, all independently.
We then have for every $e\in E(G)$ that the expected number of triangles in $\cT'$ containing $e$ is $\frac{1}{4}p^2 n$. Using Lemma~\ref{lem:chernoff}\ref{chernoff eps}, it is then easy to see that with high probability, $\cT'$ is $(n^{-1/3},p/2)$-regular.\COMMENT{\begin{align*}
\prob{|\cT'(e)|\neq (1\pm n^{-1/3}) \frac{1}{4}p^2 n} \le 2\eul^{-\frac{n^{-2/3}p^2 n}{12}} \le \eul^{-n^{1/4}}.
\end{align*}}

It remains to show that $\psi$ exists. Note that since $G$ is $(\xi,4,p)$-typical, we have $|\cT^{(3)}(e)|=(1\pm \xi)p^2 n$ and $|\cT^{(5)}(e)| = (1\pm \xi)^3 p^9 n^3 /6$ for all $e\in E(G)$. Thus, for all $e\in E(G)$, 
\begin{align}
\mbox{defining } c_e:=\frac{p^2 n-|\cT^{(3)}(e)|}{4\cdot |\cT^{(5)}(e)|} \mbox {, we have } |c_e| \le \frac{6\xi p^2 n}{4(1-\xi)^3 p^9 n^3} \le \frac{3 \xi}{ p^{7}n^{2}} .  \label{weights}
\end{align}

For every $e\in E(G)$ and $J\in \cT^{(5)}(e)$, we define a function $\psi_{e,J} \colon \cT^{(3)}\to \bR$ as follows: for $T\in \cT^{(3)}$ with $T\In J$, let
\begin{align}
\psi_{e,J}(T) := \begin{cases} 1/3, & |V(T)\cap e|\in \Set{0,2}; \\  -1/6, & |V(T)\cap e|=1; \end{cases} \label{shifter def}
\end{align}
and let $\psi_{e,J}(T):=0$ if $T\not\In J$.
Observe that for all $e'\in E(G)$, 
\begin{align}
\sum_{T \in \cT^{(3)}(e')}\psi_{e,J}(T)=\begin{cases} 1, & e'=e,\\ 0, & e'\neq e. \end{cases} \label{fractional identity}
\end{align}

We now define $\psi\colon \cT^{(3)}\to [0,1]$ as $$\psi:=\frac{1}{4} + \sum_{e\in E(G)} c_e \sum_{J\in \cT^{(5)}(e)}\psi_{e,J}.$$
For every $e\in E(G)$, we have
\begin{align*}
\sum_{T\in \cT^{(3)}(e)}\psi(T) &= \frac{1}{4}|\cT^{(3)}(e)| + \sum_{e'\in E(G)}c_{e'}\sum_{J\in \cT^{(5)}(e')}\sum_{T\in \cT^{(3)}(e)}\psi_{e',J}(T)\\
                             &\overset{\eqref{fractional identity}}{=} \frac{1}{4}|\cT^{(3)}(e)| + c_e|\cT^{(5)}(e)| \overset{\eqref{weights}}{=}\frac{1}{4}p^2 n,
\end{align*}
as desired.
Moreover, for every $T\in \cT^{(3)}$, there are at most $\binom{n}{2}\binom{5}{2}\le 5 n^2$ pairs $(e,J)$ for which $e\in E(G)$, $J\in \cT^{(5)}(e)$ and $T\In J$. Hence,
\begin{eqnarray*}
|\psi(T)-1/4| & \le &\sum_{e\in E(G),J\in \cT^{(5)}(e)\colon T\In J}|c_e| |\psi_{e,J}(T)| \overset{\eqref{weights},\eqref{shifter def}}{\le} 5n^2 \cdot \frac{3\xi}{p^7 n^2} \cdot \frac{1}{3} \le 1/4,
\end{eqnarray*}
implying that $0\le \psi(T)\le 1$ for all $T\in \cT^{(3)}$, as needed.
\endproof

\section{Concluding remarks: general designs}

Instead of working with the typicality notion as in the proof of Theorem~\ref{thm:qr}, one can also impose more specific conditions on $G$ that allow the proof to work. For instance, a natural condition is `regularity', meaning that every edge is contained in roughly the same number of triangles. Another requirement arises from the regularity  boosting step, for which we would also require that every edge lies in $\Omega(n^3)$
copies of $K_5$. 
These properties can also be formulated for $r$-uniform hypergraphs, which led to the definition of `supercomplexes' in \cite{GKLO:16}. This allowed for the treatment of the minimum degree setting (which corresponds to Theorem~\ref{thm:main}) and the quasirandom setting (which corresponds to Theorem~\ref{thm:qr}) within a unified framework.

The proof in \cite{GKLO:16} proceeds by induction on the uniformity
$r$, both for the construction of exclusive absorbers as well as in the Cover down step.
In fact, we implicitly already used induction in the 
proof of Theorem~\ref{thm:main}:
firstly, in  the cover down step for triangles, we assumed the existence of a perfect matching in an $n$-vertex graph of minimum
degree at least $n/2$ (note that a perfect matching can be viewed as a $(2,1,1)$-design). Secondly, in the simplified construction of absorbers
mentioned after the proof of Lemma~\ref{lem:transformers},
we choose a perfect matching in the link graph of a given vertex.
In turn, to prove e.g.~the existence of decompositions of hypercliques into tetrahedra
(which corresponds to the existence of $(4,3,1)$-designs), 
a strengthening of Theorem~\ref{thm:qr}
(as well as the existence of suitable perfect matchings) is used both in the construction of absorbers as well as in the Cover down step.

%
%

\begin{thebibliography}{10}

\bibitem{AKS:97}
N.~Alon, J.-H.~Kim, and J.~Spencer, \emph{Nearly perfect matchings in regular
  simple hypergraphs}, Israel J. Math.~\textbf{100} (1997), 171--187.

\bibitem{AY:05}
N.~Alon and R.~Yuster, \emph{On a hypergraph matching problem}, Graphs
  Combin.~\textbf{21} (2005), 377--384.

\bibitem{BKLMO:17}
B.~Barber, D.~K\"uhn, A.~Lo, R.~Montgomery, and D.~Osthus, \emph{Fractional
  clique decompositions of dense graphs and hypergraphs}, J. Combin. Theory
  Ser.~B~\textbf{127} (2017), 148--186.

\bibitem{BKLO:16}
B.~Barber, D.~K\"uhn, A.~Lo, and D.~Osthus, \emph{Edge-decompositions of graphs
  with high minimum degree}, Adv. Math.~\textbf{288} (2016), 337--385.

\bibitem{BKLOT:17}
B.~Barber, D.~K\"uhn, A.~Lo, D.~Osthus, and A.~Taylor, \emph{Clique
  decompositions of multipartite graphs and completion of {L}atin squares}, J.
  Combin. Theory Ser.~A~\textbf{151} (2017), 146--201.

\bibitem{DP:19}
M.~Delcourt and L.~Postle, \emph{Progress towards {N}ash-{W}illiams' conjecture
  on triangle decompositions}, arXiv:1909.00514 (2019).

\bibitem{DT:97}
D.~Dor and M.~Tarsi, \emph{Graph decomposition is {NP}-complete: a complete
  proof of {H}olyer's conjecture}, SIAM J. Comput.~\textbf{26} (1997),
  1166--1187.

\bibitem{dross:16}
F.~Dross, \emph{Fractional triangle decompositions in graphs with large minimum
  degree}, SIAM J. Discrete Math.~\textbf{30} (2016), 36--42.

\bibitem{DLR:95}
R.~A.~Duke, H.~Lefmann, and V.~R{\"o}dl, \emph{A fast approximation algorithm
  for computing the frequencies of subgraphs in a given graph}, SIAM J.
  Comput.~\textbf{24} (1995), 598--620.

\bibitem{DH:19}
P.~J.~Dukes and D.~Horsley, \emph{On the minimum degree required for a triangle
  decomposition}, arXiv:1908.11076 (2019).

\bibitem{EGP:91}
P.~Erd\H{o}s, A.~Gy\'arf\'as, and L.~Pyber, \emph{Vertex coverings by
  monochromatic cycles and trees}, J. Combin. Theory Ser. B~\textbf{51} (1991),
  90--95.

\bibitem{GKLMO:19}
S.~Glock, D.~K\"uhn, A.~Lo, R.~Montgomery, and D.~Osthus, \emph{On the
  decomposition threshold of a given graph}, J. Combin. Theory Ser.
  B~\textbf{139} (2019), 47--127.

\bibitem{GKLO:16}
S.~Glock, D.~K\"uhn, A.~Lo, and D.~Osthus, \emph{The existence of designs via iterative absorption: hypergraph {$F$}-designs for
  arbitrary~{$F$}}, Mem. Amer. Math. Soc. (to appear).

\bibitem{HR:01}
P.~E.~Haxell and V.~R\"odl, \emph{Integer and fractional packings in dense
  graphs}, Combinatorica~\textbf{21} (2001), 13--38.

\bibitem{JLR:00}
S.~Janson, T.~{\L}uczak, and A.~Ruci\'{n}ski, \emph{Random graphs},
  Wiley-Intersci. Ser. Discrete Math. Optim., Wiley-Interscience, 2000.

\bibitem{JKKO:ta}
F.~Joos, J.~Kim, D.~K\"uhn, and D.~Osthus, \emph{Optimal packings of bounded
  degree trees}, J. Eur. Math. Soc. (to appear).

\bibitem{keevash:14}
P.~Keevash, \emph{The existence of designs}, arXiv:1401.3665 (2014).

\bibitem{keevash:18}
\bysame, \emph{Counting designs}, J. Eur. Math. Soc.~\textbf{20} (2018),
  903--927.

\bibitem{keevash:18b}
\bysame, \emph{The existence of designs~{II}}, arXiv:1802.05900 (2018).

\bibitem{keevash:18c}
\bysame, \emph{Hypergraph matchings and designs}, Proc. Int. Cong. of
  Math.~\textbf{3} (2018), 3099--3122.

\bibitem{kirkman:47}
T.~P.~Kirkman, \emph{On a problem in combinatorics}, Cambridge Dublin Math.
  J.~\textbf{2} (1847), 191--204.

\bibitem{KKO:15}
F.~Knox, D.~K\"uhn, and D.~Osthus, \emph{Edge-disjoint {H}amilton cycles in
  random graphs}, Random Structures Algorithms~\textbf{46} (2015), 397--445.

\bibitem{krivelevich:97}
M.~Krivelevich, \emph{Triangle factors in random graphs}, Combin. Probab.
  Comput.~\textbf{6} (1997), 337--347.

\bibitem{KO:09}
D.~K\"uhn and D.~Osthus, \emph{The minimum degree threshold for perfect graph
  packings}, Combinatorica~\textbf{29} (2009), 65--107.

\bibitem{KO:13}
\bysame, \emph{Hamilton decompositions of regular expanders: {A} proof of
  {K}elly's conjecture for large tournaments}, Adv. Math.~\textbf{237} (2013),
  62--146.

\bibitem{montgomery:17}
R.~Montgomery, \emph{Fractional clique decompositions of dense partite graphs},
  Combin. Probab. Comput.~\textbf{26} (2017), 911--943.

\bibitem{nash-williams:70}
C.~St. J.~A.~Nash-Williams, \emph{An unsolved problem concerning decomposition
  of graphs into triangles}, In:~Combinatorial {T}heory and its {A}pplications
  {III} (P.~Erd\H{o}s, A.~R\'enyi, and V.T.~S\'os, eds.), North Holland, 1970,
  pp.~1179--1183.

\bibitem{raman:90}
R.~Raman, \emph{The power of collision: {R}andomized parallel algorithms for
  chaining and integer sorting}, In:~Foundations of software technology and
  theoretical computer science (K.~V.~Nori and C.~E.~Veni~Madhavan, eds.),
  Lecture Notes in Comput. Sci. 472, Springer, 1990, pp.~161--175.

\bibitem{rodl:85}
V.~R\"odl, \emph{On a packing and covering problem}, European J.
  Combin.~\textbf{6} (1985), 69--78.

\bibitem{RRS:06}
V.~R\"odl, A.~Ruci\'nski, and E.~Szemer\'edi, \emph{A {D}irac-type theorem for
  {$3$}-uniform hypergraphs}, Combin. Probab. Comput.~\textbf{15} (2006),
  229--251.

\bibitem{taylor:19}
A.~Taylor, \emph{On the exact decomposition threshold for even cycles}, J.
  Graph Theory~\textbf{90} (2019), 231--266.

\end{thebibliography}
%
\providecommand{\bysame}{\leavevmode\hbox to3em{\hrulefill}\thinspace}
\providecommand{\MR}{\relax\ifhmode\unskip\space\fi MR }
\providecommand{\MRhref}[2]{%
  \href{http://www.ams.org/mathscinet-getitem?mr=#1}{#2}
}
\providecommand{\href}[2]{#2}

\vspace{1.5cm}

{\footnotesize \obeylines \parindent=0pt

Ben Barber\textsuperscript{\textasteriskcentered}, Stefan Glock\textsuperscript{\S}, Daniela K\"{u}hn\textsuperscript{\textdagger}, Allan Lo\textsuperscript{\textdagger}, Richard Montgomery\textsuperscript{\textdagger}, Deryk Osthus\textsuperscript{\textdagger}

\vspace{0.3cm}
\textsuperscript{\textasteriskcentered}School of Mathematics, University of Bristol and Heilbronn Institute for Mathematical Research, Bristol, UK.

\vspace{0.3cm}
\textsuperscript{\S}Institute for Theoretical Studies, ETH Z\"urich, 8092 Z\"urich, Switzerland

\vspace{0.3cm}
\textsuperscript{\textdagger}School of Mathematics, University of Birmingham, Edgbaston, Birmingham, B15 2TT, UK

}
\vspace{0.3cm}
\begin{flushleft}
{\it{E-mail addresses}:}
\tt{b.a.barber@bristol.ac.uk, stefan.glock@eth-its.ethz.ch, [d.kuhn,s.a.lo,r.h.montgomery,d.osthus]@bham.ac.uk}
\end{flushleft}

\end{document}